\documentclass[10pt,reqno]{amsart}

\title{Integral criteria for transportation-cost inequalities}
\author{Nathael Gozlan}
\date{\today}

\usepackage{amssymb,amsmath, amsfonts, latexsym}
\usepackage[english]{babel}
\usepackage{hyperref}

\def\marginpar#1{\ignorespaces}

\newtheorem{thm}[equation]{Theorem}
\newtheorem{lem}[equation]{Lemma}
\newtheorem{prop}[equation]{Proposition}

\newtheorem{rem}[equation]{Remark}
\newtheorem*{Hyp}{Assumptions}

\numberwithin{equation}{section}
\newcommand{\st}{\text{ such that }}

\newcommand{\la}{\left\{}
\newcommand{\ra}{\right\}}
\newcommand{\lp}{\left(}
\newcommand{\rp}{\right)}
\newcommand{\lc}{\left[}
\newcommand{\rc}{\right]}

\newcommand{\ep}{\varepsilon}
\renewcommand{\d}{\delta}

\newcommand{\R}{\mathbb{R}}

\def\AArm{\fam0 \rm}%
\newdimen\AAdi%
\newbox\AAbo%
\def\AAk#1#2{\setbox\AAbo=\hbox{#2}\AAdi=\wd\AAbo\kern#1\AAdi{}}%
\newcommand{\1}{{\ensuremath{{\AArm 1\AAk{-.8}{I}I}}}}

\renewcommand{\P}{\mathbb{P}}
\newcommand{\E}{\mathbb{E}}

\newcommand{\X}{\mathcal{X}}
\newcommand{\IX}{\int_\mathcal{X}}
\newcommand{\PX}{\mathcal{P}(\X)}

\newcommand{\Hnm}{\operatorname{H}(\nu\mid\mu)}

\newcommand{\Tnm}{\mathcal{T}_c(\nu,\mu)}
\newcommand{\Tdnm}{\mathcal{T}_d(\nu,\mu)}
\newcommand{\Tsqnm}{\mathcal{T}_{d^2}(\nu,\mu)}
\newcommand{\Tdpnm}{\mathcal{T}_{d^p}(\nu,\mu)}
\newcommand{\normnm}{\|\nu-\mu\|^*_{\Phi}}
\newcommand{\chinmTV}{\|\chi\cdot(\nu-\mu)\|_{TV}}

\renewcommand{\a}{\alpha}
\newcommand{\as}{\a^\circledast}

\newcommand{\T}{\mathbb{T}}

\newcommand{\C}{\mathcal{C}}

\newcommand{\Bl}{\mathrm{BLip}_1(\X,d)}

\renewcommand{\L}{\mathbb{L}}
\newcommand{\ta}{\tau_\a}
\newcommand{\Lea}{\L_{\ta}(\X,\mu)}

\newcommand{\tatilde}{\tau_{\widetilde{\a}}}
\newcommand{\Leatilde}{\L_{\tatilde}(\X,\mu)}

\newcommand{\dom}{\mathrm{dom\ }}

\begin{document}


\address{Modal-X, Universit\'e Paris 10. B\^at. G, 200 av. de la R\'epublique. 92001 Nanterre Cedex, France}

\email{nathael.gozlan@u-paris10.fr}


\begin{abstract}
In this paper, we provide a characterization of a large class of transportation-cost inequalities in terms of exponential integrability of the cost function under the reference probability measure. Our results completely extend the previous works by Djellout, Guilin and Wu \cite{DGW03} and Bolley and Villani \cite{BV03}.
\end{abstract}

\maketitle

\section{Introduction}
In all the paper, $(\X,d)$ will be a polish space equipped with its
Borel $\sigma$-field. The set of probability measures on $\X$ will
be denoted by $\PX$.\medskip

\subsection{Norm-entropy inequalities and transportation cost inequalities}The aim of
this paper is to give necessary and sufficient conditions for inequalities of the following form :
\begin{equation}\label{norm-entropy}
\forall \nu\in \PX,\quad \a\lp\normnm\rp\leq\Hnm,
\end{equation}
where
\begin{itemize}
\item $\a : \R^+\to\R^+\cup\{+\infty\}$ is a convex lower semi-continuous (l.s.c) function vanishing at $0$,
\item The semi-norm $\normnm$ is defined by
\begin{equation}
\normnm := \sup_{\varphi\in \Phi} \la\IX \varphi\,d\nu-\IX \varphi\,d\mu\ra,
\end{equation}
where $\Phi$ is a set of bounded measurable functions on $\X$ which is symmetric, i.e.
\begin{center}$\varphi \in \Phi \Rightarrow - \varphi \in \Phi$,\end{center}
\item The quantity $\Hnm$ is the
relative entropy of $\nu$ with respect to $\mu$ defined by $$\Hnm =
\IX \log \frac{d\nu}{d\mu}\,d\nu,$$ if $\nu$ is absolutely
continuous with respect to $\mu$ and $+\infty$ otherwise.
\end{itemize}
Inequalities of the form (\ref{norm-entropy}) were introduced by C. L\'eonard and the
author in \cite{GozLeo}. They are called \emph{norm-entropy inequalities}. An important
particular case, is when $\Phi$ is the set of all bounded $1$-Lipschitz functions on $\X$ :
$\Phi=\Bl$. Indeed, in that case $\normnm$ is the optimal transportation cost between $\nu$ and $\mu$ associated to the
metric cost function $d(x,y)$. Let us recall that
if $c:\X\times\X\to\R^+$ is a lower
semi-continuous function, then the optimal transportation cost between $\nu\in
\PX$ and $\mu \in \PX$ is defined by
\begin{equation}\label{Tnm}
\Tnm = \inf\ \int_{\X^2}c(x,y)\,d\pi(x,y)
\end{equation}
where $\pi$ describes the set $\Pi(\nu,\mu)$ of all probability
measures on $\X\times\X$ having $\nu$ for first marginal and $\mu$
for second marginal. According to Kantorovich-Rubinstein duality theorem (see e.g Theorem 1.3 of \cite{Vill}),
if the cost function $c$ is the metric $d$, the following identity holds
\begin{equation}
\Tdnm = \sup_{\varphi\in \Bl} \la\IX \varphi\,d\nu-\IX \varphi\,d\mu\ra.
\end{equation}
In this setting, inequality (\ref{norm-entropy}) becomes
\begin{equation}\label{TCI-d}
\forall \nu\in \PX,\quad \a\lp\Tdnm\rp\leq\Hnm
\end{equation}
Such an inequality is called a \emph{convex transportation-cost inequality} (convex T.C.I).\medskip

\subsection{Applications of transportation-cost inequalities}After the seminal works of
K. Marton \cite{Mar86,Mar96} and M.
Talagrand \cite{Tal96a}, new efforts have been made in order to
understand this kind of inequalities. The reason of this interest is
the link between T.C.I and concentration of measure
inequalities. Namely, according to a general argument du to K.
Marton, if $\mu$ satisfies (\ref{TCI-d}), then $\mu$ has the following
concentration property
$$\forall A\subset \X \text{ s.t. }\mu(A)\geq \frac 1 2,\quad \forall
\ep\geq r,\quad \mu(A^\ep)\geq 1-e^{-\a(\ep-r)},$$
with $r=\a^{-1}(\log(2))$ and $A^\ep=\{x\in \X : d(x,A)\leq \ep\}$.
For a proof of this fact, see e.g. Theorem 9 of \cite{GozLeo}. Other
applications of T.C.Is were investigated in \cite{DGW03},
\cite{BV03}, \cite{BGV05} and \cite{GozLeo}. In these papers, it was
shown that T.C.Is are an efficient way for deriving precise
deviations results for Markov chains and empirical
processes. One can also consult \cite{CaGoz} and \cite{Goz} for applications
of norm-entropy inequalities to the study of conditional principles of
Gibbs type for empirical measures and random weighted measures.\medskip

\subsection{Necessary and sufficient conditions for norm-entropy inequalities}Our main
result gives necessary and sufficient conditions on $\mu$ for (\ref{norm-entropy}) to be
satisfied.
Before to state it, let us introduce some notations. In all what follows, $\C$ will denote
the set of convex functions $\a :\R^+\to \R^+\cup\{+\infty\}$ which are
lower semi continuous (l.s.c) and such that $\a(0)=0$. For a given $\alpha$, the monotone
convex conjugate of $\a$ will be denoted by $\as$. It is defined by
$$\forall s\geq 0,\quad \as(s)=\sup_{t\geq 0}\la st - \a(t)\ra.$$
Note that, if $\a$ belongs to $\C$, then $\as$ also belongs to $\C$. Furthermore, one has
the relation $\a^{\circledast\,\circledast}=\a$.
If $\a$ is in $\C$, the Orlicz space
$\Lea$ associated to the function $\ta:=e^\a-1$ is defined by
$$\Lea=\la f : \X\to\R\st \exists \lambda>0,\  \IX \ta\lp\frac f \lambda\rp d\mu<+\infty\ra,$$
where $\mu$ almost everywhere equal functions are identified. The
space $\Lea$ is equipped with its classical Luxemburg norm
$\|\,.\,\|_{\ta}$, i.e
$$\forall f \in \Lea,\quad \|f\|_{\ta}=\inf \la \lambda>0 \st  \IX \ta\lp\frac f \lambda\rp d\mu\leq 1 \ra.$$
We will need the following assumptions on $\a$ :
\begin{Hyp}~
\begin{itemize}
\item[$(A_1)$ :] The effective domain of $\as$ is open on the left, i.e $\la s\in \R^+ : \as(s)<+\infty\ra=[0,b[$, for some $b>0$.
\item[$(A_2)$ :] The function $\as$ is super-quadratic near $0$, i.e
\begin{equation}\label{super-quad}
\exists s_{\as}>0,c_{\as}>0,\quad
\forall s \in [0,s_{\as}],\quad \as(s)\geq c_{\as} s^2.
\end{equation}
\end{itemize}
\end{Hyp}
We can now state the main result of this paper, which will be proved in section 2.
\begin{thm} \label{norm-entropy-equiv}
Let $\alpha \in \C$ satisfy assumptions $(A_1)$ and $(A_2)$ and $\mu\in \PX$.
The following statements are equivalent :
\begin{enumerate}
\item $\displaystyle{\exists a>0 \st ,\quad\forall \nu\in \PX,\quad  \a\lp\frac \normnm a \rp\leq \Hnm}$
\item $\exists M>0\st ,\quad\forall \varphi \in \Phi,\quad \|\varphi -\langle\varphi,\mu\rangle\|_{\ta}\leq M$.
\end{enumerate}
More precisely, if (1) holds true then one can take $M=3a$.
Conversely, if (2) holds true, then one can take $a=\sqrt{2}m_\a M$, with $m_\a$ defined by
\begin{center}
$m_\a=e\max\left(\frac{1}{\a^{-1}(2)\sqrt{c_{\as}(1-u)}},\frac{1}{u}\right)$
with $u\in[0,1[$ such that : $\displaystyle{\frac{u}{\sqrt{1-u}}\leq
s_{\as}\sqrt{c_\as}\quad\text{ and } \quad\frac{u^3}{1-u}\leq
2},$\end{center}
where the constants $s_{\as}$ and $c_{\as}$ are
given by (\ref{super-quad}).
\end{thm}
\begin{rem}~
\begin{itemize}
\item If $\Phi$ contains an element which is not $\mu$-a.e constant, and if
inequality (\ref{norm-entropy}) holds for some $\a\in \C$, then $\a$ satisfies
assumption $A_2$ (see Lemma \ref{about-A_2}).
\item The constant $a=\sqrt{2}m_\a M$ is not optimal. This can be easily checked
by considering the celebrated Pinsker inequality, i.e
\begin{equation}\label{Pinsker}
\forall \nu\in \PX,\quad \frac{\|\nu- \mu\|^2_{TV}}{2}\leq \Hnm,
\end{equation}
where $\|\nu-\mu\|_{TV}$ is the total-variation norm which is defined by
$$\|\nu-\mu\|_{TV}=\sup\la\IX \varphi\,d\nu-\IX\varphi\,d\mu, |\varphi|\leq 1\ra.$$
\end{itemize}
\end{rem}
In order to prove Theorem \ref{norm-entropy-equiv}, we will take advantage of the dual
formulation of norm-entropy inequalities developed in \cite{GozLeo}. Namely, according
to Theorem 3.15 of \cite{GozLeo}, we have the following result :
\begin{thm}\label{norm-entropy-dual}
The inequality $$\forall \nu \in \PX,\quad \a\lp\frac{\normnm}{a}\rp\leq \Hnm,$$ with $\a \in \C$
is equivalent to the following condition :
\begin{equation}\label{Laplace-upper-bound}
\forall \varphi\in \Phi,\quad \forall s\in \R^+,\quad \IX e^{s\varphi}\,d\mu\leq e^{s\langle\varphi,\mu\rangle+\as(as)}.
\end{equation}
\end{thm}
According to (\ref{Laplace-upper-bound}), the only thing to know is how to majorize the
Laplace transform of a centered random variable $X$ knowing that this random variable
satisfies an Orlicz integrability condition of the form : $\E\lc e^{\a\lp\frac X \lambda\rp}\rc<+\infty$,
for some $\lambda>0$.
Estimates of this kind are very useful in probability theory, because they enable us to control
the deviation probabilities of sums of independent and identically distributed random variables.
In \cite{GozLeo}, we have shown how to deduce Pinsker inequality from the classical Hoeffding estimate
(see Section 2.3 of \cite{GozLeo}). We also proved that the weighted version of Pinsker inequality (\ref{WCPK-BV03.1}) recently obtained by
Bolley and Villani in \cite{BV03} is a consequence of Bernstein estimate (see Corollaries 3.23 and 3.24 of \cite{GozLeo}). Here,
Theorem \ref{norm-entropy-equiv} will follow very easily from the following theorem which is du
 to Kozachenko and Ostrovskii (see \cite{KO85} and \cite{BulKo} p. 63-68) :
\begin{thm}\label{KO}
Suppose that $\a \in \C$ satisfies Assumptions $(A_1)$ and $(A_2)$, then for all $f\in \Lea$
 such that $\IX f\,d\mu=0$, the following holds
$$\forall s\geq 0,\quad \IX e^{sf}\,d\mu\leq e^{\as\lp a s\rp},$$
with $a=\sqrt{2}m_\a \|f\|_{\ta}$, where $m_\a$ is the constant defined in
Theorem \ref{norm-entropy-equiv}.
\end{thm}
For further informations on the preceding result, we refer to Chapter VII of \cite{GozPhD} (p. 193-197) where a complete detailed proof is given.
Before proving Theorem \ref{norm-entropy-equiv}, we discuss below some of its applications.
\subsection{Applications to T.C.Is}
Applying the preceding theorem to the case where $\Phi$ is the Lipschitz ball $\Bl$, one obtains the following result.
\begin{thm}\label{TCI-equiv-d}
Let $\alpha \in \C$ satisfy assumptions $(A_1)$ and $(A_2)$ and $\mu\in \PX$ be such that $\IX d(x_0,x)\,d\mu(x)<+\infty$ for all $x_0\in \X$. The following statements are equivalent :
\begin{enumerate}
\item $\displaystyle{\exists a>0 \st  \forall \nu\in \PX,\quad \a\lp\frac{\Tdnm}{a}\rp\leq
\Hnm}$.
\item For all $x_0\in \X$, the function $d(x_0,\,.\,)\in \Lea$.
\end{enumerate}
More precisely, if $(2)$ holds true, then one can take $a=2\sqrt{2}m_\a
\inf_{x_0\in\X}\|d(x_0,\,.\,)\|_{\ta}$, where $m_\a$ was defined in Theorem \ref{norm-entropy-equiv}.
\end{thm}
Actually, other transportation cost inequalities can be deduced from Theorem \ref{norm-entropy-equiv}.
 Using a majorization technique developed by F. Bolley and C. Villani in \cite{BV03}, we will prove
 the following result :
\begin{thm}\label{TCI-equiv-c}Let $c(\,.\,,\,.\,)$ be a cost function such that $c(x,y)=q(d(x,y))$,
where $q:\R^+\to \R^+$ is an increasing convex function satisfying the $\Delta_2$-condition, i.e
\begin{equation}\label{Delta2}
\exists K>0,\quad \forall x \in \R^+,\quad q(2x)\leq Kq(x).
\end{equation}
If $\a \in \C$ satisfies assumptions $(A_1)$ and $(A_2)$, then for all $\mu\in \PX$ such that $\IX c(x_0,x)\,d\mu(x)<+\infty$ for all $x_0\in \X$, the following statements are equivalent :
\begin{enumerate}
\item $\displaystyle{\exists a>0,\quad \forall \nu \in \PX,\quad \a\left(\frac{\Tnm}{a}\right)\leq \Hnm,}$
\item For all $x_0\in \X$, the function $c(x_0,\,.\,)\in \Lea$.
\end{enumerate}
More precisely, if $(2)$ holds true then one can take $a=\sqrt{2}K m_\a
\inf_{x_0\in\X}\|c(x_0,\,.\,)\|_{\ta}$. Furthermore, if $\dom \a=\R^+$ then the following inequality holds
\begin{equation}\label{Main-Inequality}
\forall \nu\in \PX,\quad \Tnm \leq \sqrt{2}Km_\a\inf_{x_0\in \X,\, \d>0}\frac{1}{\d}\lp1+\frac{\log\IX e^{\d\a\lp c(x_0,x)\rp}\,d\mu(x)}{\log 2}\rp\a^{-1}\lp\Hnm\rp
\end{equation}
\end{thm}
Contrary to what happens in the case where $c$ is the metric $d$,
a transportation-cost inequality $\a\lp\Tnm\rp\leq\Hnm$ can hold even
if $\a$ does not satisfy Assumption $(A_2)$. The most known example is
Talagrand inequality, also called $\T_2$-inequality. Let us recall that
a probability measure $\mu$ on $\R^n$ satisfies the Talagrand inequality $\T_2(a)$ if
\begin{equation}\label{T2}
\forall \nu \in \PX,\quad \Tsqnm\leq a\Hnm,
\end{equation}
where $d(x,y)=\sqrt{\sum_{i=1}^n (x_i-y_i)^2}$. Gaussian measures do
satisfy a $\T_2$-inequality. This was first shown by Talagrand in
\cite{Tal96a}. In this case, the corresponding $\a$ is a linear
function and hence its monotone conjugate $\as$ does not satisfy
$(A_2)$. Sufficient conditions are known for Talagrand inequality.
In \cite{OVill00}, it was shown by F. Otto and C. Villani that if
$d\mu=e^{-\Phi}dx$ is a probability measure on $\R^n$ satisfying a
logarithmic Sobolev inequality with constant $a$, then it also
satisfies the inequality $\T_2(a)$. Furthermore, if $\mu$ satisfies
$\T_2(a)$, then it satisfies the Poincar\'e inequality with a
constant $a/2$. An alternative proof of these facts was proposed in
\cite{BGL01} by S.G. Bobkov, I. Gentil and M. Ledoux. In a recent
paper P. Cattiaux and A. Guillin gave an example of a probability
measure satisfying $\T_2$ but not the logarithmic Sobolev inequality
(see \cite{CaGui}). A necessary and sufficient condition for $\T_2$
is not yet known. Other examples of transportation-cost inequalities
involving a linear $\a$ can be found in \cite{BGL01}, \cite{GGM} and
\cite{CaGui}. The common feature of these $\T_2$-like inequalities
is that they enjoy a dimension free tensorization property (see e.g
Theorem 4.12 of \cite{GozLeo}) which in turn implies a dimension
free concentration phenomenon.
\subsection{About the literature}
Theorems \ref{TCI-equiv-c} and \ref{TCI-equiv-d} extend previous results obtained by
 H. Djellout, A. Guillin and L. Wu in \cite{DGW03} and by F. Bolley and C. Villani in \cite{BV03}.

In \cite{DGW03}, H. Djellout, A. Guillin and L. Wu obtained the first integral criteria
for the so called $\T_1$-inequality. Let us recall that a probability measure $\mu$ on $\X$
 is said to satisfy the inequality $\T_1(a)$ if
\begin{equation}\label{T1}
\forall \nu \in \PX,\quad \Tdnm^2\leq a\Hnm.
\end{equation}
According to Jensen inequality, $\Tdnm^2\leq \Tsqnm$, and thus
$\T_2(a)\Rightarrow \T_1(a)$. The inequality
$\T_1$ is weaker than $\T_2$ and it is also considerably easier to
study. According to Theorem 3.1 of \cite{DGW03}, the following propositions are equivalent :
\begin{enumerate}
\item $\exists a>0,\st  \mu \text{ satisfies } \T_1(a)$
\item $\exists \d>0 \st
\displaystyle{\int_{\X^2}e^{\d d(x,y)^2}\,d\mu(x)d\mu(y)<+\infty}$
\end{enumerate}
More precisely, if
$\displaystyle{\int_{\X^2}e^{\d d(x,y)^2}\,d\mu(x)d\mu(y)<+\infty}$
for some $\d>0$, then one can take
\begin{equation}\label{DGW2}
a=\frac{4}{\d^2}\sup_{k\geq1}\lp\frac{(k!)^2}{(2k!)}\rp^{1/k}\lc\int_{\X^2}e^{\d^2
d(x,y)^2}\,d\mu(x)d\mu(y)\rc^{1/k}<+\infty.
\end{equation}
The link between the constants $a$ and $\delta$ was then improved by
F. Bolley and C. Villani in \cite{BV03} (see (\ref{BV03-2})
bellow).\medskip

In \cite{BV03}, F. Bolley and C. Villani obtained the following weighted versions of Pinsker inequality : if $\chi:\X\rightarrow\R^+$, is a measurable function, then for all $\nu \in \PX$,
\begin{align}
\chinmTV&\leq \lp\frac 3 2 +\log\IX e^{2\chi}\,d\mu\rp\lp\sqrt{\Hnm}+\frac{1}{2}\Hnm\rp\label{WCPK-BV03.1}\\
\chinmTV&\leq \sqrt{1+\log \IX  e^{\chi^2}\,d\mu}\sqrt{2\Hnm}\label{WCPK-BV03.2}
\end{align}
Using the following upper bound (see \cite{Vill}, prop. 7.10)
\begin{equation}\label{sandwich-BV}
\Tdpnm\leq 2^{p-1}\|d(x_0,\,.\,)^p\cdot(\nu-\mu)\|_{TV},
\end{equation}
they deduce from (\ref{WCPK-BV03.1}) and (\ref{WCPK-BV03.2}) the following transportation cost inequalities involving cost functions of the form $c(x,y)=d(x,y)^p$ with $p\geq 1$ : $\forall \nu \in \PX,$
\begin{align}
\Tdpnm^{1/p}&\leq 2\inf_{x_0\in\X,\,\d>0}\lc\frac{1}{\d}\lp\frac{3}{2}+\log\IX
e^{\d d(x_0,x)^p}d\mu(x)\rp\rc^{1/p}\cdot\lc\Hnm^{1/p}+\lp\frac{\Hnm}{2}\rp^{1/2p}\rc,\label{BV03-1}\\
\Tdpnm&\leq 2\inf_{x_0\in\X,\,\d>0}\lc\frac{1}{2\d}\lp1+\log\IX e^{\d
d(x_0,x)^{2p}}d\mu(x)\rp\rc^{1/2p}\cdot \Hnm^{1/2p}.\label{BV03-2}
\end{align}
Note that for $p=1$, the constant in (\ref{BV03-2}) is sharper than (\ref{DGW2}). Note also that, up to numerical factors, (\ref{BV03-1}) and (\ref{BV03-2}) are particular cases of (\ref{Main-Inequality}).\\

In order to derive T.C.Is from norm-entropy inequalities, we will follow the lines of \cite{BV03}. To do this, we will deduce from Theorem \ref{norm-entropy-equiv} a general version of weighted Pinsker inequality (see Theorem \ref{WCPK}). Theorem \ref{TCI-equiv-c} will follow from Theorem \ref{WCPK} and from Lemma \ref{sandwich} which generalizes inequality (\ref{sandwich-BV}).
\medskip
\section{Necessary and sufficient conditions for norm-entropy inequalities.}
Let us begin with a remark on Assumption $(A_2)$.
\begin{lem}\label{about-A_2}
Suppose that $\Phi$ contains a function $\varphi_0$ which is not $\mu$-almost everywhere constant.
If $\mu$ satisfies the inequality
$$\forall \nu\in \PX,\quad \a\lp\normnm\rp\leq \Hnm,$$
then $\a$ satisfies Assumption $(A_2)$.
\end{lem}
\proof
Let us define $\Lambda_{\varphi_0}(s)=\log \IX e^{s\varphi_0}\,d\mu$, for all $s\in \R$.
According to Theorem \ref{norm-entropy-dual}, we have
$$\forall s\geq0,\quad \Lambda_{\varphi_0}(s)-s\langle\varphi_0,\mu\rangle \leq \as(s).$$
It is well known that
$$\lim_{s \rightarrow 0^+}\frac{\Lambda_{\varphi_0}(s)-s\langle\varphi_0,\mu\rangle }{s^2}=
\frac{1}{2}\operatorname{Var}_\mu(\varphi_0)>0.$$
From this follows that
$\displaystyle{\liminf_{s \rightarrow
0^+}\frac{\as(s)}{s^2}>0},$ which easily implies (\ref{super-quad}).
\endproof
\begin{rem}
Note that if all the elements of $\Phi$ are $\mu$-almost everywhere constant, then
$\normnm=0$ for all $\nu \ll \mu$. Inequality (\ref{norm-entropy}) is thus satisfied,
 for all $\a\in \C$.
\end{rem}
The rest of this section is devoted to the proof of Theorem \ref{norm-entropy-equiv}.
The following lemma will be useful in the sequel :
\begin{lem}\label{DZ}
Let $X$ be a random variable such that $\E\left[e^{\delta|X|}\right]<+\infty$, for some $\delta>0$.
Let us denote by $\Lambda_X$ the Log-Laplace of $X$, which is defined by $\Lambda_X(s)=\log\E\lc e^{sX}\rc$,
 and by $\Lambda_X^*$ its Cram\'er transform defined by $\Lambda_X^*(t)=\sup_{s\in \R}\la st - \Lambda_X(s)\ra$,
  then the following upper-bound holds :
$$\forall \ep\in [0,1[,\quad\E\lc e^{\ep \Lambda_X^*(X)}\rc \leq \frac{1+\ep}{1-\ep}.$$
\end{lem}
\proof (See also Lemma 5.1.14 of \cite{DZ}.)
Let $a<b$ with $a\in \R\cup\{-\infty\}$ and $b\in \R\cup\{+\infty\}$ be the endpoints of $\dom\Lambda_X^*$.
Since $\Lambda_X^*$ is convex l.s.c, $\{\Lambda_X^*\leq t\}$ is an interval with endpoints
$a \leq a(t)\leq b(t) \leq b$, for all $t\geq0$. As a consequence,
$$\forall t\geq0,\quad \P(\Lambda_X^*(X)>t)=\P(X<a(t))+\P(X>b(t)).$$
Let $m=\E[X]$. Since $\Lambda_X^*(m)=0$, $a(t)\leq m$. But for all $u\leq m$,
it is well known that
\begin{equation}\label{DZ1}
\P(X\leq u)\leq \exp(-\Lambda_X^*(u))
\end{equation}
If $a(t)>a$, the continuity of $\Lambda_X^*$
on $]a,b[$ easily implies that $\Lambda_X^*(a(t))=t$. Thus, according to (\ref{DZ1}),
$$\P(X<a(t))\leq e^{-t}.$$
If $a(t)=a$, then
$$\P(X<a)=\lim_{n \rightarrow +\infty}\P(X<a-1/n)
\overset{(i)}{\leq}\lim_{n \rightarrow +\infty}\exp(-\Lambda_X^*(a-1/n))
\overset{(ii)}{=}\lim_{n \rightarrow +\infty}0=0,$$
where (i) comes from (\ref{DZ1}) and (ii) from $a-1/n \notin \dom\Lambda_X^*$.\\
Therefore, in all cases $\P(X<a(t))\leq e^{-t}$. In the same way, we have
$\P(X>b(t))\leq e^{-t}.$
As a consequence,
\begin{equation}\label{DZ2}
\forall t \geq 0,\quad\P\lp\Lambda_X^*(X)>t\rp\leq 2e^{-t}.
\end{equation}
Finally, integrating by parts and using (\ref{DZ2}) in ($\ast$) bellow, we get
\begin{align*}
\E\lc e^{\ep \Lambda_X^*(X)}\rc
& = \int_{-\infty}^{+\infty}e^t\P\lp\Lambda_X^*(X)>t/\ep\rp\,dt =  \int_{-\infty}^0e^t\,dt+\int_{0}^{+\infty}e^t\P(\Lambda_X^*(X)>t/\ep)\,dt\\
&\overset{(\ast)}{\leq} 1+2\int_{0}^{+\infty}e^{(1-1/\ep)t}\,dt
=\frac{1+\ep}{1-\ep}.
\end{align*}
\endproof
Now, let us prove Theorem \ref{norm-entropy-equiv}.\\
\emph{Proof of Theorem \ref{norm-entropy-equiv}.}
Let us show that (1) implies (2). For $\varphi \in \Phi$, according to
Theorem \ref{norm-entropy-dual} and using the fact that $-\varphi\in \Phi$, we have
\begin{equation}\label{proof-Th.norm-entropy-equiv1}
\forall s\in \R,\quad \log \IX e^{s(\varphi-\langle\varphi,\mu\rangle}\,d\mu\leq \as(|as|).
\end{equation}
Define $\widetilde{\varphi}:=\varphi-\langle\varphi,\mu\rangle$ and
$\Lambda_{\widetilde{\varphi}}(s):=\log \IX e^{s(\varphi-\langle\varphi,\mu\rangle}\,d\mu$.
Equation (\ref{proof-Th.norm-entropy-equiv1}) immediately yields
$$\forall t\in \R,\quad \a\lp \frac{|t|}{a}\rp=\sup_{s\in \R}\la st - \as(|as|)\ra\leq \sup_{s\in \R}\la st- \Lambda_{\widetilde{\varphi}}(s)\ra=\Lambda_{\widetilde{\varphi}}^*(t).$$
According to Lemma \ref{DZ}, $\IX e^{\ep \Lambda_{\widetilde{\varphi}}^*(\widetilde{\varphi})}\,d\mu\leq \frac{1+\ep}{1-\ep}$, for all $\ep \in [0,1[$. Thus $\IX e^{\ep \a\lp\frac{\widetilde{\varphi}}{a}\rp}\,d\mu\leq \frac{1+\ep}{1-\ep}$. Since $\a\lp\frac{|\,.\,|}{a}\rp$ is convex and $\a(0)=0$, we have $\a\lp \frac{\ep|t|}{a}\rp\leq\ep\a\lp \frac{|t|}{a}\rp$. Therefore,
$\IX e^{\a\lp \frac{\ep|\widetilde{\varphi}|}{a}\rp} d\mu \leq \frac{1+\ep}{1-\ep}$.
In other words,
$$\forall \varphi\in \Phi,\quad \forall \ep\in [0,1[,\quad \IX \ta\lp \frac{\ep \widetilde{\varphi}}{a}\rp d\mu\leq \frac{2\ep}{1-\ep}.$$
It is now easy to see that $\left\|\widetilde{\varphi}\right\|_{\ta}\leq 3a$, for all $\varphi\in \Phi$.

\medskip
Now let us show that (2) implies (1).
According to Theorem \ref{KO},
$$\forall s\geq 0,\quad \IX e^{s\varphi}\,d\mu \leq e^{s\langle\varphi,\mu\rangle+\as\lp\sqrt{2}m_\a \|\varphi-\langle\varphi,\mu\rangle\|_{\ta}s\rp},$$
for all $\varphi \in \Phi$. As it is assumed that $\|\varphi-\langle\varphi,\mu\rangle\|_{\ta}\leq M$, for all $\varphi\in \Phi$, we thus have
$$\forall \varphi \in \Phi,\quad \forall s\geq 0,\quad \IX e^{s\varphi}\,d\mu \leq e^{s\langle\varphi,\mu\rangle+\as\lp as\rp},$$
with $a=\sqrt{2}m_\a M$. According to Theorem \ref{norm-entropy-dual}, this implies that $\mu$ satisfies the inequality
$$\forall \nu\in \PX,\quad \a\lp\frac \normnm a\rp\leq \Hnm.$$
~\hfill $\square$\\
\paragraph{\textbf{Example : Weighted Pinsker inequalities.}} Let $\chi : \X\to
\R^+$ be a measurable function and let $\Phi_\chi$ be the set of
bounded measurable functions $\varphi$ on $\X$ such that
$|\varphi|\leq \chi$. In this framework, it is easily seen that
$$\|\nu-\mu\|_{\Phi_\chi}^*=\chinmTV,$$
where $\|\gamma\|_{TV}$ denotes the total-variation of the signed
measure $\gamma$.
\begin{thm}\label{WCPK}
Suppose that $\IX \chi\,d\mu<+\infty$ and that $\alpha \in \C$
satisfies Assumptions $(A_1)$ and $(A_2)$, then the following
propositions are equivalent :
\begin{enumerate}
\item $\displaystyle{\exists a>0,\st \forall \nu \in \PX,\quad \alpha\lp\frac \chinmTV a\rp}\leq
\Hnm$,
\item $\chi \in \Lea$.
\end{enumerate}
More precisely, if $\chi \in \Lea$, then one can take
$a=2\sqrt{2}m_\alpha\|\chi\|_{\ta}$. Conversely, if (1) holds
true, then
$$\|\chi\|_{\ta}\leq \la
\begin{array}{ll}
3a,&\text{if }\mu \text{ has no atoms}\\
3a +\IX\chi\,d\mu\cdot\|\1\|_{\ta},&\text{otherwise}
\end{array}\right.$$
Furthermore, the Luxemburg norm $\|\chi\|_{\ta}$ can be
estimated in the following way :
\begin{itemize}
\item If $\dom \a=\R^+$, then $\displaystyle{\|\chi\|_{\ta}\leq \inf_{\d>0}\la\frac 1 \d \lp 1+ \frac{\log \IX e^{\a\lp\d\chi\rp}d\mu}{\log
2}\rp\ra}$\\
\item If $\dom \a=[0,r_\a[\text{ or }[0,r_\a]$, then $\Lea=\L_{\infty}(\X,\mu)$ and
$$r_a^{-1}\|\chi\|_\infty \leq\|\chi\|_{\ta}\leq \sup\la t>0 : \a(t)\leq \log 2\ra^{-1} \cdot \|\chi\|_\infty.$$
\end{itemize}
\end{thm}
\begin{rem}
If $\a \in \C$ satisfies Assumptions $(A_1)$ and $(A_2)$ and is such that $\dom \a=\R^+$, we have thus shown the following weighted version of Pinsker inequality :
\begin{equation}\label{WCPK-final-form}
\forall \nu \in \PX,\quad \chinmTV\leq 2\sqrt{2}m_\a\inf_{\d>0}\la\frac 1 \d \lp 1+ \frac{\log\IX e^{\a\lp\d\chi\rp}d\mu}{\log
2}\rp\ra\a^{-1}\lp\Hnm\rp
\end{equation}
Inequality (\ref{WCPK-final-form}) completely extends Bolley and Villani's results (\ref{WCPK-BV03.1}) and (\ref{WCPK-BV03.2}). The proof of Bolley and Villani is very different from ours. Roughly speaking, it relies on a direct comparison of the two integrals $\IX \chi\left|\frac{d\nu}{d\mu}-1\right|\,d\mu$ and $\IX \frac{d\nu}{d\mu}\log \frac{d\nu}{d\mu}d\mu$.
\end{rem}
\noindent \emph{Proof of Theorem \ref{WCPK}.}
According to Theorem \ref{norm-entropy-equiv}, it suffices to show that
\begin{equation}\label{Phi_chi}
2\|\chi\|_{\ta}\geq\sup_{\varphi \in \Phi_\chi}\la\|\varphi - \langle\varphi,\mu\rangle\|_{\ta}\ra\geq \la\begin{array}{ll}\|\chi\|_{\ta}&\text{if }\mu\text{ is non-atomic}\\
\|\chi\|_{\ta}-\IX\chi\,d\mu\cdot\|\1\|_{\ta}&\text{otherwise}.
\end{array}\right.
\end{equation}
Let us prove the first inequality of (\ref{Phi_chi}) : If $\varphi \in \Phi_\chi$, then $|\varphi|\leq \chi$, thus $\|\varphi - \langle\varphi,\mu\rangle\|_{\ta}\leq \|\chi\|_{\ta} + \|\langle\varphi,\mu\rangle\|_{\ta}$. Thanks to Jensen inequalty, for all $\lambda>0$, we have $\IX \ta\lp\frac{\langle\varphi,\mu\rangle}{\lambda}\rp d\mu\leq\IX \ta\lp\frac{\varphi}{\lambda}\rp d\mu.$ Thus, $\|\langle\varphi,\mu\rangle\|_{\ta}\leq\|\varphi\|_{\ta}$, which proves the desired inequality.

\medskip
\noindent Thanks to triangle inequality $\sup_{\varphi \in \Phi_\chi}\|\varphi - \langle\varphi,\mu\rangle\|_{\ta}\geq \|\chi- \langle \chi,\mu\rangle\|_{\ta}\geq \|\chi\|_{\ta}-\|\IX\chi\,d\mu\|_{\ta}=\|\chi\|_{\ta}-\IX\chi\,d\mu\cdot\|\1\|_{\ta}$.

\medskip \noindent Suppose that $\mu$ has no atoms, then $\chi\cdot\mu$ has no atoms too. As a consequence, there exists a measurable set $A \subset \X$ such that $\int_A \chi\,d\mu=\frac 1 2 \IX \chi\,d\mu$. Define $\widetilde{\chi}=\chi\1_A-\chi\1_{A^c}$. Then $\left|\widetilde{\chi}\right|=\chi$ and $\langle\widetilde{\chi},\mu\rangle =0$. Thus $\sup_{\varphi \in \Phi_\chi}\|\varphi - \langle\varphi,\mu\rangle\|_{\ta}\geq \|\widetilde{\chi}-\langle\widetilde{\chi},\mu\rangle\|_{\ta}=\|\widetilde{\chi}\|_{\ta}=\|\chi\|_{\ta}.$

\medskip Now, let us explain how to majorize the Luxemburg norms. Suppose that $\dom \a=\R^+$.
If $\|\chi\|_{\ta} \leq \frac{1}{\d}$ or if $\displaystyle{\IX e^{\a(\d\chi)}\,d\mu=+\infty }$, there is nothing to prove. Let us assume that $\|\chi\|_{\ta} \geq \frac{1}{\d}$ and that $\displaystyle{\IX e^{ \a(\d\chi)}\,d\mu<+\infty }$.
Then, denoting $\lambda=\|\chi\|_{\ta}$, we have
$$2^{\d\lambda}  \overset{(i)}{=} \lc\IX \exp\a\lp\frac{\chi}{\lambda}\rp d\mu\rc^{\d\lambda}
 \overset{(ii)}{\leq }  \IX \exp \d\lambda \a\lp\frac{\chi}{\lambda}\rp d\mu
 \overset{(iii)}{\leq }  \IX \exp  \a\lp\d\chi\rp d\mu $$
where (i) come from the definition of $\lambda=\|\chi\|_{\ta}$, (ii) from Jensen inequality and (iii) from the inequality $\a(x/M)\leq \a(x)/M$, for all $M\geq1$. Taking the $\log$ in both side of the above inequality yields $\lambda \leq \frac{1}{\d\log 2}\IX \exp  \a\lp\d\chi\rp d\mu$. Thus in any case,
$$\lambda \leq \frac 1 \d+\frac{1}{\d\log 2}\IX \exp  \a\lp\d\chi\rp d\mu,$$
for all $\d>0$, which is the desired results.

\medskip
The case where $\dom \a$ is a bounded interval is left to the reader.\hfill $\square$\\
\begin{rem}It is easy to show that when $\a(x)=x^2$, the Luxemburg norm $\|\chi\|_{\tau_{x^2}}$ can be estimated in the following way :
$$\|\chi\|_{\tau_{x^2}}\leq\inf_{\d>0}\frac{1}{\d}\sqrt{1+\frac{\log \IX e^{\d^2\chi^2}\,d\mu}{\log 2}}.$$
With this upper-bound, one obtains
\begin{equation}\label{WCPK-x^2}
\chinmTV\leq 2m_{x^2}\inf_{\d>0}\frac{1}{\d}\sqrt{1+\frac{\log\IX e^{\d^2\chi^2}\,d\mu}{\log 2}}\cdot\sqrt{2\Hnm},
\end{equation}
which differs from (\ref{WCPK-BV03.2}) only by numerical factors. The following proposition gives a way to improve the constants in the preceding inequality.
\end{rem}
\begin{prop}
For every measurable function $\chi:\X\to\R^+$, the following inequality holds
\begin{equation}\label{WCPK-x^2-improved}
\chinmTV\leq\inf_{\d>0}\frac{1}{\d}\sqrt{1+4\log \IX e^{\d^2\chi^2}\,d\mu}\cdot\sqrt{2\Hnm}.
\end{equation}
\end{prop}
\proof
First let us show that if $X$ is a real random variable such that $\E\lc e^{ X^2}\rc<+\infty$ one has the following upper bound :
\begin{equation}\label{sub-gaussian}
\forall s\geq 0,\quad\E\lc e^{s(X-\E[X])}\rc\leq e^{s^2/2}\cdot \E\lc e^{X^2}\rc^{2s^2}.
\end{equation}
Let $\widetilde{X}$ be an independent copy of $X$. According to Jensen inequality, we have $\E\lc e^{s(X-\E[X])}\rc\leq \E\lc e^{s(X-\widetilde{X})}\rc$. The random variable $X-\widetilde{X}$ is symmetric, thus $\E\lc (X-\widetilde{X})^{2k+1}\rc=0$, for all $k$. Consequently,
$$\E\lc e^{s(X-\E[X])}\rc \leq \E\lc e^{s(X-\widetilde{X})}\rc=\sum_{k=1}^{+\infty}\frac{s^{2k}\E\lc(X-\widetilde{X})^{2k}\rc}{(2k)!}
\leq\sum_{k=1}^{+\infty}\frac{s^{2k}\E\lc(X-\widetilde{X})^{2k}\rc}{2^k\cdot k!}=\E\lc e^{s^2(X-\widetilde{X})^2/2}\rc.$$
It is easily seen that $\E\lc e^{s^2(X-\widetilde{X})^2/2}\rc\leq\E\lc e^{s^2X^2}\rc^2,$ and if $s\leq 1$, $\E\lc e^{s^2X^2}\rc^2\leq \E\lc e^{X^2}\rc^{2s^2}.$ Hence,
$$\forall s\leq1,\quad\E\lc e^{s(X-\E[X])}\rc \leq\E\lc e^{X^2}\rc^{2s^2}.$$
But if $s\geq 1$, one has
$$\E\lc e^{s(X-\E[X])}\rc\leq \E\lc e^{s(X-\widetilde{X})}\rc\leq \E\lc e^{s^2/2 + (X-\widetilde{X})^2/2}\rc\leq e^{s^2/2}\cdot \E\lc e^{X^2}\rc^{2}\leq e^{s^2/2}\cdot \E\lc e^{X^2}\rc^{2s^2}.$$
So, the inequality
$$\E\lc e^{s(X-\E[X])}\rc\leq e^{s^2/2}\cdot \E\lc e^{X^2}\rc^{2s^2}$$
holds for all $s\geq0$.

\medskip
Let $\varphi$ be a bounded measurable function such that $|\varphi|\leq \chi$. Applying inequality (\ref{sub-gaussian}), one obtains immediately
$$\IX e^{s(\varphi-\langle\varphi,\mu\rangle)}\,d\mu\leq e^{s^2M^2/2},$$
with $M=\sqrt{1+4\log \IX e^{\d^2\chi^2}\,d\mu}.$ Thus, according to Theorem \ref{norm-entropy-dual}
the following norm-entropy inequality holds :
$$\chinmTV\leq\sqrt{1+4\log \IX e^{\chi^2}\,d\mu}\cdot\sqrt{2\Hnm}.$$
Replacing $\chi$ by $\d\chi$ and using homogeneity one obtains (\ref{WCPK-x^2-improved}).
\endproof
\begin{rem}
Note that (\ref{WCPK-x^2-improved}) is sharper than (\ref{WCPK-x^2}). But (\ref{WCPK-BV03.2}) is still sharper than (\ref{WCPK-x^2-improved}).
\end{rem}
\section{Applications to transportation cost inequalities.}
In this section, we will see how to derive transportation-cost inequalities from norm-entropy inequalities.
Let us begin with the proof of Theorem \ref{TCI-equiv-d}.

\medskip
\emph{Proof of Theorem \ref{TCI-equiv-d}.}
First let us show that $(1)$ implies $(2)$. According to Theorem \ref{norm-entropy-equiv}, one has $\sup_{\varphi\in \Bl}\|\varphi-\langle\varphi,\mu\rangle\|_{\ta}\leq 3a.$ In particular, using an easy approximation technique, $\|d(x_0,\,.\,)-\langle d(x_0,\,.\,),\mu\rangle\|_{\ta}\leq 3a$, and thus $d(x_0,\,.\,)\in \Lea$.

\medskip Now let us see that $(2)$ implies $(1)$. Let $x_0\in \X$ ; observe that $\Tdnm=\|\nu-\mu\|_{\Phi_{x_0}}$, with $\Phi_{x_0}=\la\varphi \in \Bl : \varphi(x_0)=0\ra$. But $\Phi_{x_0}\subset\widetilde{\Phi}_{x_0}:=\la\varphi : \forall x\in \X, |\varphi(x)|\leq d(x_0,x)\ra$. Thus, $\Tdnm\leq\|\nu-\mu\|_{\widetilde{\Phi}_{x_0}}=\|d(x_0,\,.\,)\cdot(\nu-\mu)\|_{TV}$. Applying Theorem \ref{WCPK}, one concludes that if $d(x_0,\,.\,)\in \Lea$, then the inequality $\forall \nu \in \PX,\quad \a\lp\frac{\Tdnm}{a}\rp\leq \Hnm$ holds with $a=2\sqrt{2}m_\a\|d(x_0,\,.\,)\|_{\ta}$. As this is true for all $x_0\in \X$, the same inequality holds for $a=2\sqrt{2}m_\a\inf_{x_0\in \X}\|d(x_0,\,.\,)\|_{\ta}$.\hfill $\square$\\

When the cost function is of the form $c(x,y)=q(d(x,y))$, we will  use the following result which is adapted from Proposition 7.10 of \cite{Vill} :
\begin{lem}
Let $c$ be a cost function on $\X$ of the form $c(x,y)=q(d(x,y))$, with $q:\R^+\to\R^+$ an increasing convex function. Let $x_0\in \X$ and define $\chi_{x_0}(x)=\frac{1}{2}q(2d(x,x_0))$, for all $x\in\X$. Then the following inequality holds :
\begin{equation}\label{sandwich}
\forall \nu \in \PX,\quad q\lp\Tdnm\rp\leq \Tnm \leq \|\chi_{x_0}\cdot(\nu-\mu)\|_{TV}.
\end{equation}
\end{lem}
\proof
For all $\pi \in \Pi(\nu,\mu)$, Jensen inequality yields
$\displaystyle{q\left(\int_{\X^2}d(x,y)\,d\pi(x,y)\right)\leq\int_{\X^2}q(d(x,y))\,d\pi(x,y)}$. Thus according to the definition of $\Tnm$ (see (\ref{Tnm})), one deduce immediately the first inequality in (\ref{sandwich}).
It follows from the triangle inequality and the convexity of $q$ that
$$
c(x,y)=q(d(x,y))\leq q(d(x,x_0)+d(y,y_0))
\leq \frac{1}{2}\left[q(2d(x,x_0))+q(2d(y,x_0))\right]
= \chi_{x_0}(x)+\chi_{x_0}(y).
$$
Thus $c(x,y)\leq d_{\chi_{x_0}}(x,y)$, with $d_{\chi_{x_0}}(x,y)=(\chi_{x_0}(x)+\chi_{x_0}(y))\1_{\{x\neq y\}}$ and consequently  $\Tnm\leq \mathcal{T}_{d_{\chi_{x_0}}}(\nu,\mu)$. But $\mathcal{T}_{d_{\chi_{x_0}}}(\nu,\mu)=\|\chi_{x_0}\cdot(\nu-\mu)\|_{TV}$ (see for instance, Prop. VI.7 p. 154 of \cite{GozPhD}), which proves the second part of (\ref{sandwich}).
\endproof
Using the second part of inequality (\ref{sandwich}) together with Theorem \ref{WCPK}, one immediately derives the following result which is the first half of Theorem \ref{TCI-equiv-c} :
\begin{prop}\label{TCI-equiv-c-sufficient}
Let $c$ be a cost function on $\X$ of the form $c(x,y)=q(d(x,y))$, with $q:\R^+\to\R^+$ an increasing convex function and $\a\in \C$ satisfying Assumptions $(A_1)$ and $(A_2)$. Then the following T.C.I holds
\begin{equation}\label{TCI-c}
\forall \nu\in \PX,\quad \a\lp\frac \Tnm a\rp\leq \Hnm,
\end{equation}
with $\displaystyle{a=\sqrt{2}m_\a\inf_{x_0\in \X}\|q(2d(x_0,\,.\,))\|_{\ta}}$.
Furthermore, if $q$ satisfies the $\Delta_2$-condition (\ref{Delta2}) with constant $K>0$, then one can take $\displaystyle{a=\sqrt{2}K m_\a\inf_{x_0\in \X}\|c(x_0,\,.\,)\|_{\ta}}$.
\end{prop}
\begin{rem}
If $q$ satisfies the $\Delta_2$-condition and $\dom \a=\R^+$, then $\mu$ satisfies the following T.C.I :
$$\forall \nu\in \PX,\quad \Tnm \leq \sqrt{2}Km_\a\inf_{x_0\in \X,\, \d>0}\frac{1}{\d}\lp1+\frac{\log\IX e^{\d\a\lp c(x_0,x)\rp}\,d\mu(x)}{\log 2}\rp\a^{-1}\lp\Hnm\rp $$
\end{rem}
Now, let us prove the second half of Theorem \ref{TCI-equiv-c} :
\begin{prop}\label{TCI-equiv-c-necessary}
Let $c$ be a cost function on $\X$ of the form $c(x,y)=q(d(x,y))$, with $q:\R^+\to\R^+$ an increasing convex function satisfying the $\Delta_2$-condition (\ref{Delta2}) with a constant $K>0$ and let $\a\in \C$ satisfy Assumption $(A_1)$. If $\IX c(x_0,x)\,d\mu(x)<+\infty$ for all $x_0\in \X$ and if the T.C.I (\ref{TCI-c}) holds for some $a>0$, then the function $c(x_0,\,.\,)$ belongs to $\Lea$ for all $x_0 \in \X$.
\end{prop}
\proof
According to the first part of inequality (\ref{sandwich}), $q\lp\Tdnm\rp\leq\Tnm$, thus, if (\ref{TCI-c}) holds for some $a>0$, then $\widetilde{\a}\lp\Tdnm\rp\leq\Hnm$, for all $\nu \in \PX$, where $\widetilde{\a}(x)=\a\lp\frac{q(x)}{a}\rp$. According to Theorem \ref{norm-entropy-equiv}, this implies that $\displaystyle{\sup_{\varphi \in \Bl}\|\varphi - \langle\varphi,\mu\rangle\|_{\tatilde}\leq 3}$. In particular, using an easy approximation argument, it is easy to see that $\|d(x_0,\,.\,)-\langle d(x_0,\,.\,),\mu\rangle\|_{\tatilde}\leq 3$, which implies that $d(x_0,\,.\,)\in \Leatilde$. Let $\lambda>0$ be such that $\IX \tatilde\lp\frac{d(x_0,x)}{\lambda}\rp d\mu(x)<+\infty$ and let $n$ be a positive integer such that $2^n\geq \lambda$. Then, according to the $\Delta_2$ condition satisfied by $q$, one has $q\lp \frac x \lambda \rp\geq q\lp \frac{x}{2^n}\rp\geq \frac{1}{K^n}q(x)$, for all $x \in \R^+$. Consequently, $\tatilde\lp\frac{x}{\lambda}\rp\geq \ta\lp\frac{q(x)}{aK^n}\rp$, for all $x\in \R^+$. From this follows that
$$\IX \ta\lp\frac{c(x_0,x)}{aK^n}\rp d\mu(x)=\IX \ta\lp\frac{q(d(x,x_0))}{aK^n}\rp d\mu(x)\leq\IX \tatilde\lp\frac{d(x_0,x)}{\lambda}\rp d\mu(x)<+\infty$$
and thus $c(x_0,\,.\,)\in \Lea$.
\endproof
\emph{Proof of Theorem \ref{TCI-equiv-c}.} Theorem \ref{TCI-equiv-c} follows immediately from Propositions \ref{TCI-equiv-c-sufficient} and \ref{TCI-equiv-c-necessary}.\hfill $\square$\\
\bibliographystyle{plain}
\bibliography{bib}
\end{document}